\documentclass[11pt]{amsart}
\usepackage[T1]{fontenc}
\usepackage[utf8]{inputenc}
\usepackage{lmodern,amsmath,amssymb}
\usepackage[hidelinks]{hyperref}
\newtheorem{theorem}{Theorem}[section]
\newtheorem{lemma}[theorem]{Lemma}
\newtheorem{proposition}[theorem]{Proposition}

\theoremstyle{definition}
\newtheorem{remark}[theorem]{Remark}
\newtheorem{notation}[theorem]{Notation}

\newcommand{\m}{\mathfrak m}
\newcommand{\n}{\mathfrak n}

\newcommand{\q}{\mathfrak q}
\newcommand{\D}{\mathcal D}
\newcommand{\M}{\mathcal M}
\newcommand{\N}{\mathcal N}

\DeclareMathOperator{\ann}{ann}
\DeclareMathOperator{\Ass}{Ass}
\DeclareMathOperator{\Supp}{Supp}
\DeclareMathOperator{\Spec}{Spec}
\DeclareMathOperator{\Hom}{Hom}
\DeclareMathOperator{\End}{End}
\DeclareMathOperator{\ord}{ord}

\title{On the Matlis duals of $\mathcal D$-modules in characteristic zero}
\author{Tu\u gba Mahmut\c cepo\u glu}
\address{Department of Mathematics, Gebze Technical University,
Gebze, Kocaeli, Turkey}
\email{t.mahmutcepoglu@gtu.edu.tr}

\begin{document}
\begin{abstract}
Let $R=k[[x_1,\ldots,x_n]]$ be a formal power series ring over a field $k$ of characteristic zero, and let $\D=\D(R,k)$ be the ring of $k$-linear differential operators. We prove that if $\M$ is a nonzero left $\D$-module with $(0)\notin\Ass_R(\M)$, then
$(0)\in\Ass_R D(\M)$, where $D(-)$ denotes Matlis duality, without assuming holonomicity or finite generation over $\D$. As a significant consequence of this result, we establish the Lyubeznik--Y{\i}ld{\i}r{\i}m conjecture in equicharacteristic zero. We also obtain a criterion for the Matlis dual of a $\D$-module to be non-holonomic.
\end{abstract}
\maketitle
\section{Introduction}
Throughout this paper, let $R=k[[x_1,\ldots,x_n]]$ be a formal power series ring over a field $k$ of characteristic zero, where $n\geq 1$, and let $\m=(x_1,\ldots,x_n)$ be its maximal ideal. Set $E=E_R(k)$ and write $D(M)=\Hom_R(M,E)$ for the Matlis dual of an $R$-module $M$. Let $\D=\D(R,k)$ denote the ring of $k$-linear differential operators on $R$. All $\D$-modules are considered as left $\D$-modules.

\par

For a Noetherian local ring $(A,\mathfrak n)$ and an ideal $I\subseteq A$, let $H_I^i(A)$ denote the $i$-th local cohomology module of $A$ with support in $I$, and write $D_A(-)$ for Matlis
duality over $A$. Motivated by Hellus's work on associated primes and set-theoretic complete intersections \cite{Hel07}, the author (as Y\i ld\i r\i m) and Lyubeznik in  \cite{LY} conjectured that, for every regular local ring $A$ and every nonzero ideal $I\subseteq A$, the condition $H_I^i(A)\neq 0$ implies $(0)\in\Ass_A D_A(H_I^i(A))$. In the same work, the conjecture was established in equicharacteristic $p>0$ using Frobenius methods. Further results and related formulations have been obtained in \cite{BE, Zha23}.  However, Datta et al. provided a counterexample showing that the conjecture fails in general in mixed characteristic \cite[Remark 4.7]{DSZ23}.

In this paper, we establish the Lyubeznik--Y{\i}ld{\i}r{\i}m conjecture in equicharacteristic zero. Our approach uses the natural $\D$-module structure of local cohomology modules. Our main result is the following theorem on associated primes of Matlis duals of $\D$-modules, which requires neither holonomicity nor finite generation over $\D$.
\begin{theorem}\label{t1}
Let $\M$ be a nonzero $\D$-module. If $(0)\notin\Ass _R(\M)$, then $(0)\in\Ass _R D(\M)$. 
\end{theorem}
The proof of Theorem \ref{t1} combines algebraic techniques with topological methods. After reducing to a cyclic $\D$-submodule, we use the order filtration to equip its Matlis dual with a complete metrizable linear topology. The key difficulty is to pass from the existence of elements with arbitrarily small annihilators to the existence of a single element with zero annihilator.  We overcome this problem by showing that the sets of elements whose annihilators are contained in $\m^r$ are open and dense in this topology, and then applying the Baire category theorem.

We note that the Baire category theorem has already been used in commutative algebra. Indeed, Sharp and V\'amos \cite{SV} applied it to prove countable prime avoidance in complete Noetherian local rings. In this paper, we apply it to the Matlis dual of a cyclic $\D$-module equipped with the linear topology induced by the order filtration.

Theorem \ref{t1} applies to the nonzero local cohomology module $H_{I\widehat A}^i(\widehat A)$ over the completion $\widehat A$. Using flat base change and Matlis duality, we obtain
\[
H_I^i(A)\neq 0
\quad\Longrightarrow\quad
(0)\in\Ass_A D_A(H_I^i(A))
\]
for every regular local ring $A$ containing a field of characteristic zero and every nonzero ideal $I\subseteq A$. This establishes the Lyubeznik--Y{\i}ld{\i}r{\i}m conjecture in equicharacteristic zero.

Finally, we obtain a criterion for the Matlis dual of a $\D$-module to be non-holonomic. If $\M$ is a nonzero $\D$-module such that $(0)\notin\Ass _R(\M)$ and $\operatorname{fd}_R \M\leq n-2$, then $D(\M)$ is not holonomic. This follows from Theorem \ref{t1}, the equality $\operatorname{id}_R D(\M)=\operatorname{fd}_R \M$, and the injective-dimension bound of Switala and Zhang \cite[Theorem 6.2]{SZ}.
\section{Preliminaries}
In this section, we recall basic facts about $\D$-modules and the linear topology on the Matlis dual of a cyclic $\D$-module induced by the order filtration.
\par
We first recall some basic facts about $\D$-modules, and we refer the reader to \cite{Bjork} for details. The ring of $k$-linear differential operators on $R$ is the noncommutative $k$-algebra $\D=\D(R,k)=R\langle\partial_1,\ldots,\partial_n\rangle\subseteq\End_k(R)$, where $\partial_i=\partial/\partial x_i$. We have the following relations:
\begin{enumerate}
\item $\partial_i r-r\partial_i=\partial_i(r)$ for every $r\in R$, and
 \item $\partial_i\partial_j=\partial_j\partial_i$ for all $1\leq i,j\leq n$. 
\end{enumerate}
In particular, $\partial_i x_j-x_j\partial_i=\delta_{ij}$. A left $\D$-module is an $R$-module equipped with a compatible left action of $\D$.

Every nonzero differential operator $P\in \D$ has a unique expression as a finite sum $P=\sum_{\alpha\in\mathbb N^n}r_\alpha\partial^\alpha$, where $r_\alpha\in R$ and $\partial^\alpha=\partial_1^{\alpha_1}\cdots\partial_n^{\alpha_n}$. Writing $|\alpha|=\alpha_1+\cdots+\alpha_n$, the order of $P$ is defined by $\ord(P)=\max\{|\alpha|:r_\alpha\neq0\}$. For $i\geq0$, set
\begin{equation*}
\D_i=\sum_{|\alpha|\leq i}R\partial^\alpha \quad \text {and}\quad \mathcal D_{-1}=0.
\end{equation*}
Then $\{\D_i\}_{i\geq0}$ is the order filtration on $\D$, and its associated graded ring $\operatorname{gr}\D=\bigoplus_{i\geq0}\D_i/\D_{i-1}$ is isomorphic to $R[\xi_1,\ldots,\xi_n]$, where $\xi_j$ is the principal symbol of $\partial_j$.

A nonzero finitely generated $\D$-module is called holonomic if its dimension, computed from a good filtration, is $n$; see \cite{Bjork}. For every ideal $I\subseteq R$ and every $i\geq0$, the local cohomology module $H_I^i(R)$ carries a natural left $\D$-module structure and, when nonzero, is holonomic \cite[Section~2]{Ly93}. Moreover, the Matlis dual of any left $\D$-module admits a natural left $\D$-module structure \cite[Proposition 4.17]{Sw}.

We now briefly recall some standard facts about the topology defined by a filtration, which we will apply to the Matlis dual of a cyclic $\D$-module using the order filtration.

Let $M$ be an $R$-module and let $K_0\supseteq K_1\supseteq\cdots$ be a descending sequence of submodules of $M$. We say that a subset $U\subseteq M$ is \emph{open} if every $x\in U$ admits an $i\ge 0$ such that $x+K_i\subseteq U$. This defines a topology on $M$, which is called the \emph{linear topology} determined by $(K_i)_{i\ge 0}$, where the cosets $x+K_i$ form a basis of open sets. Each $K_i$ is both open and closed, and this topology is Hausdorff if and only if $\bigcap_{i\ge 0}K_i=0$; see \cite[\S 8]{MatCRT}.

Assume now that $\bigcap_{i\ge 0}K_i=0$, and let $K_{-1}=M$. For $z\in M$
we define
\[
  \nu(z)=\sup\{\,i\ge -1 : z\in K_i\,\},
\]
so that $\nu(z)=\infty$ if and only if $z=0$. Since each $K_i$ is a submodule of $M$, we have $\nu(-z)=\nu(z)$ and $\nu(z+w)\ge\min\{\nu(z),\nu(w)\}$ for all $z,w\in M$. Now define
\[
  d(x,y)=2^{-\nu(x-y)}
\]
for $x,y\in M$, where $2^{-\infty}=0$. Then $d$ is a metric on $M$ such that $d(x,z)\le\max\{d(x,y),\,d(y,z)\}$ for all $x,y,z\in M$. Moreover, we have $d(x,y)\le 2^{-i}$ if and only if $x-y\in K_i$, for each $i\ge 0$. Hence the closed balls of $d$ are the cosets of the submodules $K_i$, and so $d$
induces the linear topology determined by $(K_i)_{i\ge 0}$.

If, in addition, the canonical map $M\to\varprojlim_i M/K_i$ is an isomorphism, then every Cauchy sequence in $M$ converges, and so $(M,d)$ is a complete metric space; see \cite[pp.~102--105]{AM}.

We will use the Baire category theorem: every countable intersection of open dense subsets of a complete metric space is dense \cite[Theorem~48.2]{Munkres}. In particular, such an intersection is nonempty whenever the space is nonempty.

\begin{notation}\label{n}
Let $\N$ be a nonzero cyclic left $\D$-module. Then $\N=\D u$ for some $0\neq u\in\N$. Now consider the order filtration of $\D$ and write
\[
  F_i=\D_iu=\sum_{|\alpha|\le i}R\,\partial^{\alpha}u
\]
for each $i\ge 0$. Clearly $\N=\bigcup_{i\ge 0}F_i$ and
\[
  F_0\subseteq F_1\subseteq F_2\subseteq\cdots.
\]
Furthermore, each $F_i$ is finitely generated as an $R$-module. Now define
$\rho_i\colon D(\N)\longrightarrow D(F_i)$, for each $i\ge 0$, to be
the natural restriction homomorphism, and write $K_i=\ker\rho_i$. Then
\[
  K_0\supseteq K_1\supseteq K_2\supseteq\cdots.
\]
In the rest of this paper, we consider $D(\N)$ equipped with the linear topology determined by the descending sequence $(K_i)_{i\ge 0}$, where the cosets $\varphi+K_i$ form a basis of open sets.
\end{notation}
We now collect some well-known properties of this topology, and we include a proof for the convenience of the reader.
\begin{proposition}\label{p}
Let the notations be as in Notation~\ref{n}. Then the following hold:
\begin{enumerate}
\item[(a)] $\rho_i$ is surjective and $D(\N)/K_i\cong D(F_i)$, for
           each $i\ge 0$.
\item[(b)] $\bigcap_{i\ge 0}K_i=0$.
\item[(c)] The canonical map gives an isomorphism
           \[
             D(\N)\cong\varprojlim_{i}D(\N)/K_i.
           \]
\end{enumerate}
In particular, $D(\N)$ is a complete metric space with respect to the metric $d$ defined above. Thus the Baire category theorem can be applied to $D(\N)$.
\end{proposition}
\begin{proof}
(a) Since $E$ is injective, every $R$-linear map $F_i\to E$ extends to $\N$, so $\rho_i$ is surjective. Hence it follows from the first isomorphism theorem that $D(\N)/K_i\cong D(F_i)$.

(b) Since $\N=\bigcup_{i\ge 0}F_i$, any
$\varphi\in\bigcap_{i\ge 0}K_i$ vanishes on every $F_i$, and hence on
$\N$. Therefore, $\bigcap_{i\ge 0}K_i=0$.

(c) Applying the Matlis dual functor to the filtration $(F_i)_{i\ge 0}$, we
obtain the canonical isomorphism
\[
  D(\N)\cong\varprojlim_{i}D(F_i),
\]
where $D(F_i)\longrightarrow D(F_{i-1})$ is the natural surjection induced by
the inclusion $F_{i-1}\hookrightarrow F_i$, for each $i\ge 1$. Indeed, a compatible family of $R$-linear maps $F_i\to E$ defines a unique $R$-linear map $\N\to E$. Then the assertion follows from (a).

Finally, by (b) and (c), the discussion preceding Notation \ref{n} shows that $D(\N)$ is a complete metric space with respect to $d$.
\end{proof}

We note that, in this topology, a subset $U\subseteq D(\N)$ is open if and only if every $\varphi\in U$ admits an $i$ such that $\varphi+K_i\subseteq U$. It is dense if and only if every coset
$\varphi+K_i$ meets $U$.

\section{Proof of Main Result}
In this section, we first establish the auxiliary results on annihilators and $\D$-modules needed for the proof of Theorem \ref{t1} and then give a proof of Theorem \ref{t1}.

We start with the following well-known lemma:

\begin{lemma}[Chevalley's theorem {\cite[Exercise 8.7]{MatCRT}}]\label{chevalley}
Let $(A,\n)$ be a complete Noetherian local ring, and let $J_1\supseteq J_2\supseteq\cdots$ be a descending sequence of ideals with $\bigcap_{t\geq1}J_t=0$. For every $r\geq1$, there exists $t_0\geq1$ such that $J_t\subseteq\n^r$ for all $t\geq t_0$. Equivalently, the linear topology on $A$ defined by $(J_t)_{t\geq 1}$ is finer than or equal to the $\n$-adic topology.
\end{lemma}
The following lemma concerns the annihilator ideals under localization.

\begin{lemma}\label{l1}
Let $A$ be a Noetherian ring, $N$ be an $A$-module, and $\q\in\Spec A$. Then, for every $v\in N_{\mathfrak q}$, there exists $u\in N$ such that
\[
\ann_A(u)
=
\ann_{A_{\q}}(v)\cap A.
\]
\end{lemma}
\begin{proof}
Let $v=w/s$, where $w\in N$ and $s\notin\mathfrak q$, and set $J=\ann_{A_{\mathfrak q}}(v)\cap A$. Since $J$ is finitely generated and every element of $J$ annihilates
$w/1$, there exists $a\notin\q$ such that $aJw=0$. Let now $u=aw$. Then $J\subseteq\ann_A(u)$. Conversely, if $ru=0$, localizing at $\q$ and using that $a$ and $s$ are units gives $rv=0$, so $r\in J$. Hence, we obtain $\ann_A(u)=J$, as desired.
\end{proof}


The following lemma establishes a containment property for annihilator ideals with respect to the $\mathfrak m$-adic filtration.
\begin{lemma}\label{l2}
Let $M$ be an $R$-module, $a,b\in M$, and $r\ge1$. If
$\ann_R(b)\subseteq\m^{2r-1}$, then either $ \ann_R(a)\subseteq\m^r$ or $ \ann_R(a+b)\subseteq\m^r.$
\end{lemma}

\begin{proof}
Suppose on the contrary that neither $ \ann_R(a)\subseteq\m^r$ nor $ \ann_R(a+b)\subseteq\m^r.$ Then there exist $c\in\ann_R(a)\setminus\m^r$ and $d\in\ann_R(a+b)\setminus\m^r$. Since $ca=0$ and $d(a+b)=0$, we have $cd(a+b)=cdb=0$. Thus $cd\in\ann_R(b)\subseteq\m^{2r-1}$.

Now define $o(f)=\max\{t\geq0:f\in\m^t\}$, for a nonzero element $f\in R$. Since $R=k[[x_1,...,x_n]]$, $\operatorname{gr}_{\m}(R)\cong k[X_1,\ldots,X_n]$ is an integral domain and so $o(-)$ is additive on the products of nonzero elements. As $c,d\notin\m^r$, we obtain
\[
o(cd)
=o(c)+o(d)
\leq 2r-2.
\]
Hence we have  $cd\notin\m^{2r-1}$ which contradicts the above inclusion.
\end{proof}
We now state and prove the following lemma which shows that, for any $\D$-module $\M$, the condition $(0)\notin\Ass_R(\M)$ forces $\m \M=\M$.

\begin{lemma}\label{l3}
Let $\M$ be a left $\D$-module such that $(0)\notin\Ass_R(\M)$. Then $\m\M=\M$.
\end{lemma}

\begin{proof}
Let $0\neq u\in \M$. Since $(0)\notin\Ass_R(\M)$, there exists a nonzero element $f\in R$ such that $fu=0$. As $u\neq0$, the element $f$ cannot be a unit, so $f\in\m$. Let $d=o(f)$, where $o(-)$ is as defined in Lemma \ref{l2}, and choose a multi-index $\alpha\in\mathbb N^n$ with $|\alpha|=d$ such that the coefficient of $x^\alpha$ in $f$ is nonzero. Since $\operatorname{char}(k)=0$, the element $\partial^\alpha(f)$ has a nonzero constant term and is therefore a unit in $R$. Moreover, for every $\beta\leq\alpha$ with $\beta\neq\alpha$, we have $|\beta|<d$, so $\partial^\beta(f)\in\mathfrak m$. Applying $\partial^\alpha$ to $fu=0$ and using the Leibniz rule, we obtain
\[
0=\partial^\alpha(fu)
=\partial^\alpha(f)u+
\sum_{\substack{\beta\leq\alpha\\\beta\neq\alpha}}
\binom{\alpha}{\beta}
\partial^\beta(f)\partial^{\alpha-\beta}(u)
\]
and so we have 
\[
\partial^\alpha(f)u=-
\sum_{\substack{\beta\leq\alpha\\\beta\neq\alpha}}
\binom{\alpha}{\beta}
\partial^\beta(f)\partial^{\alpha-\beta}(u).
\]
Since every term in the sum belongs to $\m\M$, it follows that $\partial^\alpha(f)u\in\m\M$. As $\partial^\alpha(f)$ is a unit, we conclude that $u\in\m \M$. Hence $\M\subseteq\m\M$. Since the reverse inclusion is trivial, we have $\m\M=\M$, as desired.
\end{proof}

We now turn to annihilator ideals of elements in $\D$-modules. Using the local injectivity
properties of these modules, we obtain a sequence of primary annihilator ideals with zero intersection.

\begin{proposition}\label{p1} Let $\M$ be a nonzero left $\D$-module and let $0\neq\q\in\Ass_R(\M)$. Then there exist elements $\lambda_t\in \M$ and
a descending sequence of $\q$-primary ideals $J_t$ such that
\[
\ann_R(\lambda_t)=J_t
\quad\text{and}\quad
\bigcap_{t\geq1}J_t=0.
\]
In particular, for every $r\geq1$, there exists $\lambda\in \M$ such that
$\ann_R(\lambda)\subseteq\m^r$.
\end{proposition}

\begin{proof}
Let $A=R_{\q}$, $\mathfrak a=\q A$, and $T=(H_{\q}^0(\M))_{\q}$. Since $\q\in\Ass_R(\M)$, there exists $u\in \M$ with $\ann_R(u)=\q$. Thus $u/1$ is a nonzero element of $T$ annihilated by $\mathfrak a$, giving an embedding $\kappa(\q)=\operatorname{Frac}(R/{\q R})\hookrightarrow T$. By \cite[Lemma 3.4]{SZ}, $T$ is an injective $A$-module. Therefore the inclusion $\kappa(\q)\hookrightarrow T$  extends to a homomorphism $E_A(\kappa(\q))\to T$, which is injective since $\kappa(\q)$ is essential in its injective hull. Hence we may regard $E_A(\kappa(\q))$ as a submodule of $T$.

Let now $h=\operatorname{ht}(\mathfrak q)$ and choose a regular system of parameters $z_1,\ldots,z_h$ of $A$ and let $B_t=(z_1^t,\ldots,z_h^t)A$, for each $t\geq1$. Since $A/B_t$ is an Artinian Gorenstein local ring, its socle is one-dimensional over $\kappa(\q)$, that is $\operatorname{Soc}(A/B_t)\cong \kappa(\q)$, and is essential in $A/B_t$. Therefore an embedding $\kappa(\q)\hookrightarrow E_A(\kappa(\q))$ extends to an injective homomorphism $A/B_t\to E_A(\kappa(\q))$. Therefore we may consider  $A/B_t$  as a submodule of $T$ and so there exists $v_t\in T$ such that
\[
\ann_A(v_t)=B_t.
\]

Then it follows from Lemma \ref{l1} that there exists $\lambda_t\in \M$ such that $\ann_R(\lambda_t)=B_t\cap R$, for each $t\geq1$. Let now $J_t=B_t\cap R$. Since $B_t$ is $\mathfrak a$-primary, its contraction $J_t$ is $\q$-primary. Moreover, the ideals $\{B_t\}_{t\geq 1}$, and hence
$\{J_t\}_{t\geq 1}$, form descending sequences. As $B_t\subseteq\mathfrak a^t$, it follows from Krull's intersection theorem that $\bigcap_{t\geq1}B_t=0$. Since $R$ is a domain, the natural map $R\to A$ is injective, and therefore $\bigcap_{t\geq1}J_t=0$. Finally, by Lemma \ref{chevalley}, for every $r\geq1$, there exists $t\geq1$ such that $J_t\subseteq\mathfrak m^r$, as desired.
\end{proof}
\begin{remark}
Proposition \ref{p1} implies that for every nonzero left $\D$-module $\M$ and every $r\geq1$, there exists $u\in \M$ such that $\ann_R(u)\subseteq\m^r$. Indeed, if $(0)\in\Ass_R(\M)$, then $\M$ contains an element with zero annihilator. If, on the other hand, $(0)\notin\Ass_R(\M)$, the assertion follows by choosing a nonzero associated prime of $\M$ and applying Proposition \ref{p1}. However, the existence of elements with arbitrarily small annihilators does not necessarily imply that $(0)$ is an associated prime. For instance, every element of $E$ is killed by some power of  $\m$ but $\Ass_R(E)=\{\m\}$. 
\end{remark}
\begin{lemma}\label{l4}
Let the notations be as in Notation \ref{n}  and assume that $(0)\notin\Ass_R(\N)$. Then, for every $i\geq 0$ and every $r\geq 1$, there exists $b\in K_i$ such that $\ann_R(b)\subseteq\m^r.$
\end{lemma}

\begin{proof}
 Since $(0)\notin\Ass_R(\N)$, it follows from Lemma \ref{l3} that $\m\N=\N$. Then the adjointness of Hom and tensor yields the following isomorphism
\[
(0:_{D(\N)}\m)
\cong
\operatorname{Hom}_R(\N/\m\N,E)=0,
\]
from which we conclude that $\m\notin\Ass_R(D(\N)).$

Now choose $\q\in\Ass_R(D(\N))$. If $\q=(0)$, there exists $\lambda\in D(\N)$ with
$\ann_R(\lambda)=0$. If, on the other hand, $0\neq\q\in\Ass_R(D(\N))$, then it follows from Proposition \ref{p1} that there exists $\lambda\in D(\N)$ such that $J=\ann_R(\lambda)$ is $\q$-primary and $J\subseteq\m^r.$ Thus, in either case, $J$ is $\q$-primary and $J\subseteq\m^r$.

Let $y\in\m\setminus\q$. Since $F_i$ is finitely generated over $R$, the image of $\rho_i(\lambda)$ has finite length and so is killed by a power of $\m$, that is, $\mathfrak m^t\rho_i(\lambda)=0$ for some $t>0.$ 

Finally, let $b=y^t\lambda$. Then since $\rho_i(b)=\rho_i(y^t\lambda)=y^t\rho_i(\lambda)=0,$ $b\in \ker \rho_i= K_i$. Moreover, $$\ann_R(b) =(J:y^t)=J\subseteq\m^r.$$ Hence we conclude that $b$ is the desired element in $D(\N)$.
\end{proof}
\begin{lemma}\label{l5} Let the notations be as in Notation \ref{n}. Then, for any $\varphi \in D(\N)$ and $r\geq 1$, if $\ann_R(\varphi)\subseteq\m^r,$  there exists $i_0\geq 0$ such that $\ann_R(\rho_i(\varphi))\subseteq\m^r$ for every $i\geq i_0$.
\end{lemma}
\begin{proof}
Define $J=\ann_R(\varphi)$ and $J_i=\ann_R(\rho_i(\varphi))$. Since $\N=\bigcup_{i\geq 0}F_i$, the ideals $J_i$ form a descending sequence satisfying $\bigcap_{i\geq 0}J_i=J.$ Then,
 by applying Lemma \ref{chevalley} to the descending sequence $J_i/J$ in the complete local ring $R/J$, we obtain $J_i\subseteq J+\m^r =\m^r$ for all sufficiently large $i$.
\end{proof}
\begin{lemma}\label{l6} Let the notations be as in Notation \ref{n} and assume that $(0)\notin\Ass_R(\N)$. Then the subset
\[U_r=\{ \varphi\in D(\N): \ann_R(\varphi)\subseteq\m^r \}
\]
is open and dense in $D(\N)$ for every $r\geq 1$.
\end{lemma}
\begin{proof}
Let $\varphi\in U_r$. Then it follows from  Lemma \ref{l5} that $\ann_R(\rho_i(\varphi))
\subseteq\m^r$ for some $i\geq 0$. But then, for every $\psi\in\varphi+K_i$, we have $\rho_i(\psi)=\rho_i(\varphi).$ Therefore 
 
\[
\ann_R(\psi)\subseteq \ann_R(\rho_i(\psi))=\ann_R(\rho_i(\varphi))\subseteq\m^r.
\]
Hence $\varphi+K_i\subseteq U_r$ which proves that $U_r$ is open.

Let now $\varphi+K_i$ be an arbitrary basic open subset of $D(\N)$. Then, by Lemma \ref{l4}, there exists $b\in K_i$ such that $\ann_R(b)
\subseteq\m^{2r-1}.$ But then it follows from Lemma \ref{l2} that either $\ann_R(\varphi)
\subseteq\m^r$ or $\ann_R(\varphi+b)\subseteq\m^r.$ Since both $\varphi$ and $\varphi+b$
belong to $\varphi+K_i$, we have
\[
(\varphi+K_i)\cap U_r\neq\emptyset.
\]
Hence $U_r$ is dense in $D(\N)$, as desired.
\end{proof}
Now we are ready to give the proof of our main theorem.
\begin{proof}[Proof of Theorem \ref{t1}] Let $\M$ be a nonzero $\D$-module such that $(0)\notin \Ass_R(\M)$. Consider the cyclic submodule $\N=\D u$ of $\M$ for some $0\neq u\in\M$. 
Since $\N\subseteq \M$ and $(0)\notin \Ass_R(\M)$, we have $(0)\notin\Ass_R(\N)$.
Now consider the linear topology on $D(\N)$ introduced in Notation \ref{n}. Then it follows from Lemma \ref{l6} that the subset
\[
U_r=\{\varphi\in D(\N):
\ann_R(\varphi)\subseteq\m^r\}
\]
is open and dense in $D(\N)$, for every $r\geq 1$. But then since $D(\N)$ is a nonempty complete metric space by Proposition \ref{p}, the Baire category theorem implies that $\bigcap_{r\geq 1}U_r\neq\emptyset.$
Choose $\psi\in \bigcap_{r\geq 1}U_r$.  Then $\ann_R(\psi)\subseteq\m^r$ for every $r\geq 1$ and so $\ann_R(\psi)\subseteq\bigcap_{r\geq 1}\m^r=0$. Hence $\ann_R(\psi)=0$.

On the other hand, since $E$ is injective, the $R$-linear map $\psi:\N\to E$ extends to an $R$-linear homomorphism $\varphi:\M\to E$. As $\varphi|_\N=\psi$, we have $\ann_R(\varphi)
\subseteq\ann_R(\psi)=0.$ Therefore, $(0)\in\Ass_R( D(\M)),$ as desired.
\end{proof}

\begin{remark}
The assumption $(0)\notin\operatorname{Ass}_R(\M)$ in Theorem \ref{t1} cannot be omitted. Indeed, taking $\M=R$, we have $(0)\in\operatorname{Ass}_R(\M)$, whereas $D(\M)\cong E$ satisfies $\operatorname{Ass}_R(E)=\{\mathfrak m\}$.
\end{remark}
\section{Applications}
We begin with a significant application of Theorem~\ref{t1}, which completely resolves the Lyubeznik--Y{\i}ld{\i}r{\i}m conjecture in equicharacteristic zero.
\begin{theorem}\label{t2}
Let $(A,\n$) be a regular local ring containing a field of characteristic zero, and let $I\subseteq A$ be a nonzero ideal. If $H_I^i(A)\neq0$, then $(0)\in\operatorname{Ass}_A D_A(H_I^i(A))$. In particular, $\operatorname{Supp}_A D_A(H_I^i(A))=\operatorname{Spec}(A)$.
\end{theorem}
\begin{proof}
Let $\widehat A$ denote the $\n$-adic completion of $A$. As in \cite[Remark 4.1.1]{Hel07}, we reduce to the complete case using Matlis duality and flat base change. Indeed,
\[
D_A(H_I^i(A))
\cong D_{\widehat A}(H_{I\widehat A}^i(\widehat A))
\]
as $\widehat  A$-modules. Since $\widehat A$ is faithfully flat over $A$, we have $H_{I\widehat A}^i(\widehat A)\neq0$. Moreover, $\widehat A$ is a formal power series ring over a field of characteristic zero, and $H_{I\widehat A}^i(\widehat A)$ is $I\widehat A$-torsion. Thus
\[
(0)\notin\operatorname{Ass}_{\widehat A}
(H_{I\widehat A}^i(\widehat A)).
\]
Then, by Theorem \ref{t1}, $D_{\widehat A}(H_{I\widehat A}^i(\widehat A))$ contains an element with zero annihilator over $\widehat A$. Since $A\hookrightarrow\widehat A$, the same element has zero annihilator over $A$. Hence $(0)\in\operatorname{Ass}_A D_A(H_I^i(A)),$ as desired.
\end{proof}

\begin{remark}\label{r3} Theorem 5.7 of \cite{BE} appears to contradict our result, but its conclusion fails in the following elementary example:

Let $A=\mathbb C[x,y]$, $\n=(x,y)$, and $M=H_{(x)}^1(A) =\bigoplus_{j\ge1}\mathbb C[y]x^{-j}.$
Since $\partial_y:M\to M$ is surjective, we have
\[
H_{\mathrm{dR}}^2(M)
=M/(\partial_xM+\partial_yM)=0.
\]
Since $(x)$ is homogeneous, it follows from \cite[Lemma 1.4 and Remark 2.3]{RWZ22} that
\[
\Hom_A \bigl(A/\n,H_{\n}^0(M^*)\bigr)=0,
\]
where $M^*$ denotes the holonomic dual, namely the left $\D(A,\mathbb C)$-module obtained from the right module $\operatorname{Ext}^2_{\D(A,\mathbb C)} (M,\D(A,\mathbb C))$ by transposition. Therefore, $H_{\n}^0(M^*)=0$, since every nonzero $\n$-torsion module has nonzero socle. Thus all hypotheses of \cite[Theorem 5.7]{BE} are satisfied. However, for $\widehat A=\mathbb C[[x,y]]$, one has
\[
(0)\in\Ass_{\widehat A}
D_{\widehat A}\!\left(H_{(x)}^1(\widehat A)\right)
\]
by \cite[Example 2.4]{HS07}, contrary to its conclusion. Hence that statement does not provide a counterexample to Theorem \ref{t1}.
\end{remark}


We end this paper with a criterion for the Matlis dual of a $\D$-module to be non-holonomic. Recall that $\operatorname{fd}_R(M)$ denotes the flat dimension of $M$ over $R$, that is, the least length of a flat resolution of $M$, or $\infty$ if no finite flat resolution exists.
\begin{theorem}
Let $\M$ be a nonzero left $\D$-module such that $(0)\notin\Ass_R(\M)$. If $D(\M)$ is holonomic, then $\operatorname{fd}_R \M\geq n-1$. Consequently, if $\operatorname{fd}_R \M\leq n-2$, then $D(\M)$ is not holonomic.
\end{theorem}
\begin{proof}
By Theorem \ref{t1}, $\Supp_R D(\M)=\Spec(R)$, so $\dim\Supp_R D(\M)=n$. If $D(\M)$ is holonomic, the injective-dimension bound of Switala and Zhang \cite[ Theorem 6.2]{SZ} gives $\operatorname{id}_R D(\M)\geq n-1$. On the other hand, $\operatorname{id}_R D(\M)=\operatorname{fd}_R \M$. Indeed, for every $R$-module $N$ and every $j\geq0$, tensor--Hom adjunction and the injectivity of $E$ give
\[
\operatorname{Ext}_R^j(N,D(\M))
\cong D(\operatorname{Tor}_j^R(N,\M)),
\]
and $D(-)$ is faithful. Thus $\operatorname{fd}_R \M\geq n-1$, as required.
\end{proof}
For instance, if $n\ge 3$ and $0\neq f\in\m$, then $\operatorname{fd}_R H^1_{(f)}(R)\le 1$, and so $D(H^1_{(f)}(R))$ is not holonomic; see also \cite[Example 2.9]{SZ18} for the case $f=x_1$.
\section*{Acknowledgements}
The author used generative AI tools (ChatGPT, OpenAI; Claude, Anthropic) for assistance in locating the example in Remark 4.2, in searching the literature, and for language editing. All mathematical statements and proofs were independently verified by the author, who takes full responsibility for the content of the paper.

\end{document}